\documentclass[a4paper,11pt]{article}
\usepackage{amsfonts,amssymb,latexsym,amsmath}
\usepackage{multicol}
\begin{document}
\Large
\newcommand{\la}{\lambda}
\newcommand{\gog}{{\mathfrak g}}
\newcommand{\nog}{{\mathfrak n}}
\newcommand{\mog}{{\mathfrak m}}
\newcommand{\hog}{{\mathfrak h}}
\newcommand{\dog}{{\mathfrak d}}
\newcommand{\bog}{{\mathfrak b}}
\newcommand{\ut}{{\mathfrak u}{\mathfrak t}}
\newcommand{\gl}{{\mathfrak g}{\mathfrak l}}
\newcommand{\Dp}{\Delta_+}

\newcommand{\GL}{\mathrm{GL}}
\newcommand{\ldeg}{\mathrm{ldeg}}
\newcommand{\red}{\mathrm{red}}
\newcommand{\UT}{\mathrm{UT}}
\newcommand{\spann}{\mathrm{span}}
\newcommand{\eps}{\varepsilon}
\newcommand{\ad}{{\mathrm{ad}}}
\newcommand{\Ad}{{\mathrm{Ad}}}
\newcommand{\Ab}{{\Bbb A}}
\newcommand{\Nb}{{\Bbb N}}
\newcommand{\BC}{{\cal B}}
\newcommand{\Rb}{\Bbb R}

\newcommand{\ZC}{{\cal Z}}
\newcommand{\AC}{{\cal A}}
\newcommand{\DC}{{\cal D}}
\newcommand{\YC}{{\cal Y}}
\newcommand{\LC}{{\cal L}}
\newcommand{\De}{\Delta}

\newcommand{\BF}{B_i^{(1)}}
\newcommand{\BFF}{B_i^{(1.1)}}
\newcommand{\BFS}{B_i^{(1.2)}}
\newcommand{\BS}{B_i^{(2)}}
\newcommand{\BT}{B_i^{(3)}}
\newcommand{\BFO}{B_i^{(4)}}

\newcommand{\SC}{{\cal S}}
\newcommand{\IC}{{\cal I}}
\newcommand{\UC}{{\cal U}}
\newcommand{\FC}{{\cal F}}
\newcommand{\PC}{{\cal P}}

\date{}
\title{Diagram method in research on coadjoint orbits}
\author{A.N.Panov
\thanks{The paper is supported by RFBR grants  08-01-00151 and 06-01-00037}}
 \maketitle

\begin{abstract}

We correspond to any factor  algebra of the  unitriangular Lie
algebra with respect to a regular ideal some permutation. In terms
of this permutation one can construct a diagram, that allows to
calculate index and maximal dimension of coadjoint representation.

\end{abstract}

\section*{0. Introduction}

Coadjoint orbits play an important role in representation theory,
 symplectic geo\-met\-ry, mathematical physics. According to the orbit
 method of A.A.Kirillov ~\cite{K-Orb, K-62} for nilpotent Lie groups there exists one to one correspondence
between coadjoint orbits  and irreducible representations in Hilbert
spaces.  This gives a possibility to solve the problems of
representation theory and harmonic analysis in geometrical terms of
the orbit space. From  the other point of view,   coadjoint orbits
are symplectic manifolds and many of hamiltonian systems  of the
classical mechanics may be realized on these orbits
~\cite{G1,G2,G3}. However the problem of classification of all
coadjoint orbits for specific Lie groups (such as the group of
unitriangular matrices) is an open problem   that is still
interesting ~\cite{K-2003}. In the origin paper on the orbit method
~\cite{K-62} the description of algebra of invariants and
clas\-si\-fi\-ca\-tion of orbits of maximal dimension was obtained.

In the paper  ~\cite{P1} all coadjoint orbits for the groups
$\UT(n,K)$ of the size less or equal to seven were classified. In
the same paper we have got a classification of all subregular
coadjoint orbits for an arbitrary  $n$.

In the paper  ~\cite{P2} we consider the families of coadjoint
orbits associated with  involutions. The special case is a family of
orbits of maximal dimension that is associated with the involution
of maximal length. We obtain a formula of dimension of these orbits
and construct  generators of the defining ideal of the orbit; for
the canonical forms in these orbits we construct a polarization.

These paper is a continuation of the paper  ~\cite{P3} in which for
any Lie algebra  $\LC$, defined at the end of introduction, we
constructed the diagram  $\DC_\LC$. By this diagram one can easily
calculate the  maximal dimension of coadjoint orbit and the index of
 $\LC$.

Recall the index of a Lie algebra is a minimal dimension of the
stabilizer of a linear form on this Lie algebra. For nilpotent Lie
algebras the field of invariants of the coadjoint representation is
a pure transcendental extension of the main field with the degree
equal to the index of this Lie algebra Ли \cite{Dix}.

In ~\cite{P3} the diagram is constructed as a result of induction
procedure which is convenient in special cases and is not convenient
for proof of general theorems.  One should like to obtain some
formula that shows what symbols took place in the diagram.

 This problem is solved in this paper. To each considered Lie algebra  $\LC$ we correspond a permutation
  $w_\LC$.
 The main results are formulated in terms of decom\-po\-si\-tions of  $w_\LC$ into a product of reflections and its action
 in the segment  $[1,n]$ of positive integers (see theorems  2.2, 2.6 and 2.7).
    In the sequel of the paper we formulate a conjecture on structure of field of invariants of the coadjoint representation of  $\LC$.

Turn to presentation of the content of  paper.  Let $N=\UT(n,K)$ be
the group of unitriangular matrices  of size  $n\times n$ with units
on the diagonal and  with entries in the
 field $K$ of zero characteristic.
 The Lie algebra  $\nog=\ut(n,K)$ of this group consists of lower triangular matrices of size  $n\times n$ with zeros on the diagonal.
 One can define the natural representation of the group  $N$ in the conjugate space $\nog^*$ by the formular  $\Ad_g^*f(x)=f(\Ad_g^{-1} x)$,
 where  $f\in\nog^*$, $x\in\nog$ and $g\in N$.
 This representation is called   a coadjoint representation.
 We identify the symmetric algebra  $S(\nog)$ with the algebra of regular functions   $K[\nog^*]$ on the conjugate space  $\nog^*$.
Let us also identify $\nog^*$ with the subspace of upper triangulat
matrices with zeros on the diagonal. The pairing  $ \nog$ and
$\nog^*$  is realized due to the Killing form  $(a,b)=
\mathrm{Tr}(ab)$, where  $a\in \nog$, $ b\in
 \nog^*$. After this identification the coadjoint action  may be realized by the formula  $\Ad_g^*b=P(\Ad_g b)$, where $P$ is the natural
 projection of the space  of  $n\times n$-matrices onto  $\nog^*$.

Recall that for any Lie algebra  $\gog$ the algebra  $K[\gog^*]$ is
a Poisson algebra with respect to the Poissson bracket such that
$\{x,y\}=[x,y]$ for any  $ x,y\in\gog$. In the case $k=\Rb$ the
symplectic leaves with respect to this Poisson bracket
 coincides with the orbits of coadjoint representation ~\cite{K-Orb}. Respectively,
the algebra  of Casimir elements in  $K[\gog^*]$ coincides with the
algebra of invariants  $K[\gog^*]^N$ of the coadjoint
representation.

The coadjoint orbits of the group  $N$ are closed with respect to
Zariski topology in  $\nog^*$, since all orbits of a regular action
of an arbitrary algebraic unipotent group  in an affine algebraic
variety are closed  ~\cite[11.2.4]{Dix}.

To simplify  language we  shall give the following definition: a
root is an arbitrary pair $(i,j)$, where $i,j$ are positive integers
from 1 to $n$ and $i\ne j$. The permutation group  $S_n$ acts on the
set of roots by the formula  $w(i,j)=(w(i),w(j))$.

A root  $(i,j)$ is positive if  $ i > j$. Respectively, a root is
negative if  $i<j$. We denote the set set positive roots by $\Dp$.

For any root  $\eta=(i,j)$ we denote by  $-\eta$ the root  $(j,i)$.
We define the partial operation of addition  on the set of positive
roots: if $\eta=(i,j)\in \Dp$ and $\eta'=(j,m)\in \Dp$, then
$\eta+\eta'=(i,m)$.

Consider the standard basis  $\{y_{ij}:~ (i,j)\in \Dp\}$ in the
algebra
 $\nog$. We shall also use the notation $y_\xi$ for $y_{ij}$, where
 $\xi= (i,j)$.

Fix some subset  $M\subset \Dp$, that satisfies the following
condition: if in a sum of two positive roots one of summands belongs
to  $M$, then the sum also belongs to $M$. Denote by $\mog$ the
subspace spanned by  $\{y_{ij},~ (i,j)\in M\}$. By definition of
$M$, the subspace  $\mog$ is an ideal in the Lie algebra  $\nog$.

  Denote by  $\LC$ the factor algebra  Ли $\nog/\mog$ and by  $L$
 the corresponding factor group of  $N$ with respect to normal subgroup  $\exp(\mog)$.
Note that the conjugate space  $\LC^*$ is a subspace in $\nog^*$
which consists of   $f\in\nog^*$ that annihilate  $\mog$. The
coadjoint $L$-orbit for  $f\in\LC^*$ coincides with its $N$-orbit.

\section*{1. Method of construction of the diagram }

 As we say above in the paper  ~\cite{P3}  we correspond to the Lie
 algebra   $\LC$ a diagram $\DC_\LC$.
Let us state the method of construction of the diagram $\DC_\LC$ and
formulate the main assertions of the paper ~\cite{P3}. Consider the
order  $\succ$ on the set  $\Dp$ such that
$$(n,1)\succ
(n-1,1)\succ\ldots\succ
 (2,1)\succ(n,2)\succ\ldots\succ(3,2)\succ\ldots\succ(n,n-1).$$

 By the ideal $\mog$ we construct the diagram that is a  $n\times
n$-matrix,  in which all places  $(i,j)$, ~$i\leq j$, are not filled
and all other places (i.e. places in  $\Dp$) are filled by the the
symbols "$\otimes$"{}, ~"$\bullet$"{}, "$+$"\ and "$-$"\ according
to the following rules. The places  $(i,j)\in M$ are filled by the
symbol "$\bullet$"{}. Let us say that  this procedure is the zero
step in construction of the diagram.

 We put the symbol "$\otimes$"{} on the greatest( in the sense of order  $\succ$)  place  in $\Dp\setminus M$.
  Note that this symbol will take place in the first column, if the set of pairs of the form $(i,1)$ in $\Dp\setminus M$ is not empty.
Suppose that we put the symbol "$\otimes$"{} on the place $(k,t)$,~
$k>t$. Further, on all places  $(k,a)$, ~ $t< a <k$, we put the
symbol "$-$"\ and on all places $(b,t)$,~ $1< b <k$, we put the
symbol  "$+$". This procedure finishes the first step of
construction of diagram.

Further, we put the symbol "$\otimes$"{} on the greatest (in the
sense of order $\succ$) empty place in  $\Dp$. As above we put the
symbols "$+$"\ и "$-$"\ on the empty places taking into account the
following:  we put the symbols  "$+$"\ and "$-$"\ in pairs; if the
both places  $(k,a)$ and  $(a,t)$, where $ k>a>t$ are empty, we put
"$-$"\ on the first place and "$+$"\ on the second place; if one of
these places $(k,a)$ or $(a,t)$ are  already filled, then we do not
fill the other place. After this procedure we finish the step that
we call a second step.

Continuing the procedure further we have got  the diagram. We denote
this diagram by  $\DC_\LC$. The number of last step is
equal to the number of symbols "$\otimes$"{} in the diagram.\\
{\bf Example 1}. Let $n=7$, $\mog = Ky_{51}\oplus Ky_{61}\oplus
Ky_{71}\oplus Ky_{62}$. The corresponding diagram is as follows

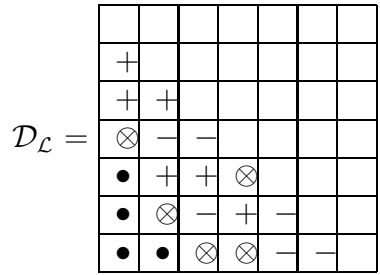
\begin{center}
{\large $\DC_\LC$ =
\begin{tabular}
{|p{0.1cm}|p{0.1cm}|p{0.1cm}|p{0.1cm}|p{0.1cm}|p{0.1cm}|p{0.1cm}|}
\hline  &  &  &  & & &  \\
\hline $+$& & &  & & &  \\
\hline $+$&$+$ & & & & &  \\
\hline $\otimes$ & $-$ & $-$  &  & &  & \\
\hline $\bullet$ & $+$ & $+$&$\otimes$ & & &\\
\hline  $\bullet$ & $\otimes$ & $-$ &$+$  &$-$ & &\\
\hline   $\bullet$&  $\bullet$ & $\otimes$ &$\otimes$  &$-$  &$-$ & \\
\hline
\end{tabular}}
\end{center}

We construct the diagram in 5 steps, beginning  with zero step:

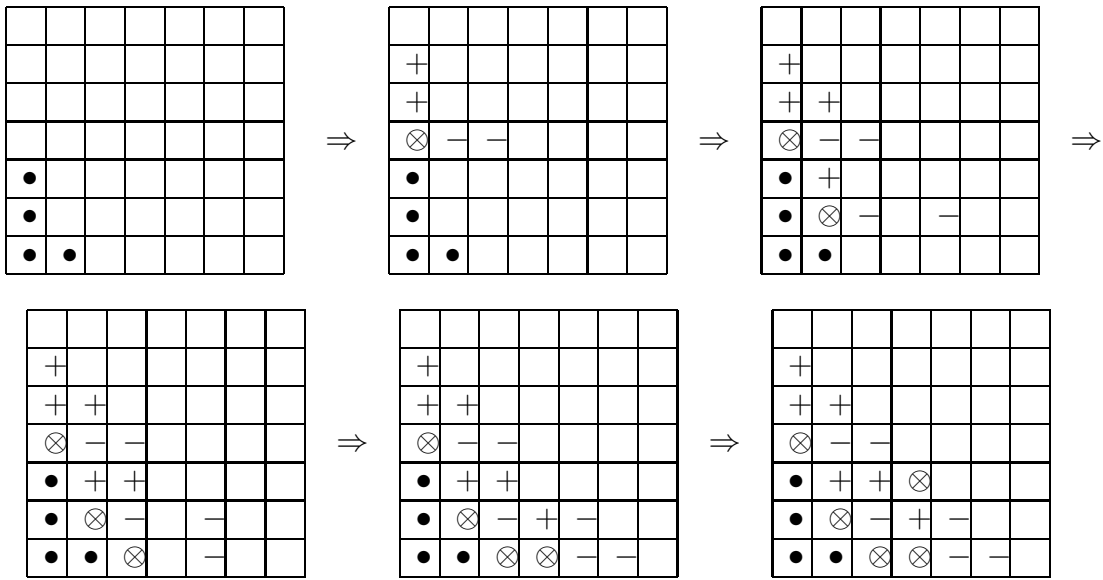
\begin{center}
{\large
\begin{tabular}{|p{0.1cm}|p{0.1cm}|p{0.1cm}|p{0.1cm}|p{0.1cm}|p{0.1cm}|p{0.1cm}|}
\hline  &  &  &  & & &  \\
\hline & & &  & & &  \\
\hline & & & & & &  \\
\hline  &  &   &  & &  & \\
\hline $\bullet$ &  & & & & &\\
\hline  $\bullet$ &  &  &  & & &\\
\hline   $\bullet$&  $\bullet$ &  &  &  & & \\
\hline
\end{tabular}
\quad $\Rightarrow$\quad
\begin{tabular}{|p{0.1cm}|p{0.1cm}|p{0.1cm}|p{0.1cm}|p{0.1cm}|p{0.1cm}|p{0.1cm}|}
\hline  &  &  &  & & &  \\
\hline $+$& & &  & & &  \\
\hline $+$& & & & & &  \\
\hline $\otimes$ & $-$ & $-$  &  & &  & \\
\hline $\bullet$ &  & & & & &\\
\hline  $\bullet$ &  &  &  & & &\\
\hline   $\bullet$&  $\bullet$ &  &  &  & & \\
\hline
\end{tabular}\quad $\Rightarrow$\quad
 \begin{tabular}{|p{0.1cm}|p{0.1cm}|p{0.1cm}|p{0.1cm}|p{0.1cm}|p{0.1cm}|p{0.1cm}|}
\hline  &  &  &  & & &  \\
\hline $+$& & &  & & &  \\
\hline $+$&$+$ & & & & &  \\
\hline $\otimes$ & $-$ & $-$  &  & &  & \\
\hline $\bullet$ & $+$ & & & & &\\
\hline  $\bullet$ & $\otimes$ & $-$ &  &$-$ & &\\
\hline   $\bullet$&  $\bullet$ &  &  &  & & \\
\hline
\end{tabular}\quad $\Rightarrow$}

\end{center}

\begin{center}
{\large
\begin{tabular}{|p{0.1cm}|p{0.1cm}|p{0.1cm}|p{0.1cm}|p{0.1cm}|p{0.1cm}|p{0.1cm}|}
\hline  &  &  &  & & &  \\
\hline $+$& & &  & & &  \\
\hline $+$&$+$ & & & & &  \\
\hline $\otimes$ & $-$ & $-$  &  & &  & \\
\hline $\bullet$ & $+$ & $+$& & & &\\
\hline  $\bullet$ & $\otimes$ & $-$ &  &$-$ & &\\
\hline   $\bullet$&  $\bullet$ & $\otimes$ &  &$-$  & & \\
\hline
\end{tabular}\quad $\Rightarrow$\quad
\begin{tabular}{|p{0.1cm}|p{0.1cm}|p{0.1cm}|p{0.1cm}|p{0.1cm}|p{0.1cm}|p{0.1cm}|}
\hline  &  &  &  & & &  \\
\hline $+$& & &  & & &  \\
\hline $+$&$+$ & & & & &  \\
\hline $\otimes$ & $-$ & $-$  &  & &  & \\
\hline $\bullet$ & $+$ & $+$& & & &\\
\hline  $\bullet$ & $\otimes$ & $-$ &$+$  &$-$ & &\\
\hline   $\bullet$&  $\bullet$ & $\otimes$ &$\otimes$  &$-$  &$-$ & \\
\hline
\end{tabular}\quad $\Rightarrow$\quad
 \begin{tabular}{|p{0.1cm}|p{0.1cm}|p{0.1cm}|p{0.1cm}|p{0.1cm}|p{0.1cm}|p{0.1cm}|}
\hline  &  &  &  & & &  \\
\hline $+$& & &  & & &  \\
\hline $+$&$+$ & & & & &  \\
\hline $\otimes$ & $-$ & $-$  &  & &  & \\
\hline $\bullet$ & $+$ & $+$&$\otimes$ & & &\\
\hline  $\bullet$ & $\otimes$ & $-$ &$+$  &$-$ & &\\
\hline   $\bullet$&  $\bullet$ & $\otimes$ &$\otimes$  &$-$  &$-$ & \\
\hline
\end{tabular}\quad\quad}
\end{center}

Denote by  $S$ (resp. $C^+$, $C^-$) the set of pairs  $(i,j)$,
filled in the diagram by the symbol  "$\otimes$" (resp. "$+$",
"$-$"). The set  $\Dp$ of positive roots  decomposes into a union of
disjoint subsets: $\Dp = M\sqcup C^+\sqcup C^-\sqcup S$.

Denote by  $\Ab_m$  the Poisson algebra $K[p_1,\ldots,p_m;
q_1,\ldots,q_m]$ with the bracket $\{p_i,q_j\}= \delta_{ij}$.

Recall that a Poisson algebra  $\AC$ is a tensor product of two
Poisson algebras  $\BC_1\otimes\BC_2$, if $\AC$ is isomorphic to
$\BC_1\otimes\BC_2$ as commutative associative algebra  and   $\{\BC_1,\BC_2\}=0$.\\
The next theorems 1.1 and 1.2 are the main results of the paper  ~\cite{P3}.\\
{\bf Theorem  1.1}~\cite{P1}. \emph{There exist  $z_1,\ldots,z_s\in
K[\LC^*]^L$, where
$s=|S|$ such that  \\
1) any  $z_i=y_{\xi_i}Q_i+P_{>i}$, where $Q_i$ is some  product of
powers of    $z_1,\ldots, z_{i-1}$ and $P_{>i}$ is a polynomial in
$\{y_\eta\}$,  $\eta\succ \xi_i$;
\\
2) denote by  $\ZC$ the set of denominators generated by
$z_1,\ldots,z_s$; the localization  $K[\LC^*]_\ZC$ of the algebra
$K[\LC^*]$ with respect to the set of denominators $\ZC$ is
isomorphic as a Poisson algebra to the tensor product
 $K[z_1^\pm,\ldots,z_s^\pm]\otimes \Ab_m$ for some  $m$.}
\\
Theorem 1.1 directly implies the following \\
{\bf Theorem  1.2}~\cite{P1}. \\ \emph{1) The field of invariants
$K(\LC^*)^L$ coincides with the field
$K(z_1,\ldots,z_s)$.\\
2) The maximal dimension of a coadjoint orbit in  $\LC^*$ equals to
the number of symbols
 "$+$"\ и "$-$"\  in the diagram  $\DC_\LC$.\\
 3) Index of Lie algebra  $\LC$ coincides with the number of symbols "$\otimes$"{} in the
 diagram.}

We have to remark that the elements  $\{z_i\}$ are constructed by
induction in the proof of theorem 1.1; this leaves unsolved the
question of finding the exact formula for these elements (for
instance, in terms of coefficients of characteristic matrix).  We
return to the question of an exact formula for generators of the
field of invariants in the sequel section of this paper.

Let us formulate two auxiliary statements from  ~\cite{P3}, that we
shall use in this paper. Denote by  $B_i$ the set of pairs $(a,b)$,~
$a>b$ that remains unfilled after the  $i$th step in the procedure
of filling the diagram. The subsets  $B_i$ form the chain:
$$B_0\supset B_1\supset\ldots\supset B_s=\emptyset,$$
where $s=|S|$.   Denote $A_i= B_i\sqcup M$,~ $\nog_i=\mathrm{span}\{
y_\eta:~\eta\in A_i\}$,~ $\LC_i=\nog_i/\mog$. Here $\Dp = A_0$.

Let  $S=\{\xi_1\succ \ldots \succ \xi_s\}$. Recall that the place
$\xi_i\in S$ is filled by the symbol "$\otimes$"\ at the $i$th step.

For $1\le i\le s$ we denote by   $C^-_i$ (resp. $C^+_i$) the subset
of pairs $(a,b)$,~ $a>b$, that are filled by the symbol "$-$"\
(resp."$+$") at the $i$th step.
\\
{\bf Proposition  1.3}~ \cite[Lemma~1]{P3}. The subspace  $\nog_i$
(resp. $\LC_i$) in $\nog$ (resp. $\LC$) is a Lie subalgebra.

For any  $1\le i\le s$ we denote
$$D^-_i = \{\eta\in \Dp|~ \xi_i \succ
\eta \quad\mbox{and}\quad \eta \in\bigsqcup_{1\le j\le i} C^-_j\}.$$

$$D^+_i = \{\eta\in \Dp|~ \xi_i \succ
\eta \quad\mbox{and}\quad \ \eta \in\bigsqcup_{1\le j\le i}
C^+_j\}.$$ Remark that the places  $\eta\in D^+_i$ are situated in
the same column as $\xi_i$.

By  $\dog^-_i$, where $1\le i\le s$,  we denote a linear subspace in
$\nog$, spanned by the vectors  $y_\eta$, where $\eta\in D^-_i$.\\
{\bf Proposition  1.4}~\cite[Lemma~2]{P3}. For any  $1\le i\le s$
the subspace  $\dog^-_i$ is a Lie subalgebra in  $\nog$.

 In the example 1 we have
$$\xi_1=(4,1),\quad \xi_2=(6,2),\quad \xi_3=(7,3), \quad
\xi_4=(7,4),\quad \xi_5=(5,4);$$
$$\begin{array}{ll}
 C^-_1=\{(4,2),(4,3)\},&\quad\quad \dog^-_1 =
 \mathrm{span}\{y_{42},y_{43}\};\\
  C^-_2=\{(6,3),(6,5)\},&\quad\quad  \dog^-_2 =
 \mathrm{span}\{y_{42},y_{43},y_{63}, y_{65}\};\\
  C^-_3=\{(7,5)\}, &\quad\quad \dog^-_3 =
 \mathrm{span}\{y_{43}, y_{63}, y_{65}, y_{75}\};\\
  C^-_4=\{(7,6)\},&\quad\quad \dog^-_4 =
 \mathrm{span}\{y_{65},y_{75},y_{76}\};\\
  C^-_5= \emptyset, &\quad\quad \dog^-_5 = \dog^-_4.\end{array}
  $$

\section*{2. Associated permutation and its decompositions }

 We correspond to a Lie algebra  $\LC$ the permutation defined as follows.\\
{\bf Definition  2.1}. Denote by $w = w_\LC$ the permutation of
$S_n$ such that \\
1) ~ $w(1)=\max\{1\le i\le n|~ (i,1)\notin M\}$; \\
2) ~ $ w(t)=\max\{1\le i\le n|~ (i,t)\notin M,
~~i\notin\{w(1),\ldots,w(t-1)\} \}$ for all  $2\le t\le n$.

As usual we denote by  $l(w)$ the number of  multipliers in
decomposition of  $w$ into products of simple  reflections. The
number  $l(w)$ coincides with the number of inversions in the
rearrangement $(w(1),\ldots,w(n))$.
\\
{\bf Theorem   2.2}. \emph{The number $l(w)$ coincides with  $\dim \LC$}.\\
{\bf Proof}. Denote by  $l(w)^{(t)}$  the number of  $k$, that $k>t$
and  $w(k)<w(t)$. Show that  $l(w)^{(t)}$ coincides with
$\dim\LC^{(t)}$, where $\LC^{(t)}=\spann \{y_{it}|~(i,t)\notin M\}$.

Let  $a^{(t)}$ be the greatest number such that the pair
$(a^{(t)},t)$ do not lie in  $M$. Obviously, $a^{(t)}\ge t$ and
$\dim\LC^{(t)}=a^{(t)}-t$.

On the other hand, the set  $\{w(1),\ldots, w(t)\}$ contains in the
segment  $[1,a^{(t)}]$. According to definition of $w$, any $c$ of
the segment  $[1,a^{(t)}]$, that do not lie in  $\{w(1),\ldots,
w(t)\}$, has  the form  $c=w(k)$, for some  $k>t$ such that
$w(k)<w(t)$. The segment  $[1,a^{(t)}]$ decomposes:
$$[1,a^{(t)}]=
\{w(1),\ldots, w(t)\}~\bigsqcup~ \{w(k)|~ k>t,~w(k)<w(t)\}.$$

Hence, $l(w)^{(t)}=a^{(t)}-t=\dim \LC^{(t)}$. Since
$$l(w)=\sum_{t=1}^n l(w)^{(t)}\quad \mbox{and}\quad  \LC=\oplus_{t=1}^n
\LC^{(t)},$$  then  $l(w)=\dim \LC$. $\Box$

 Recall that  $S=\{\xi_1\succ \xi_2 \succ\ldots \succ \xi_s\}$.
Put   $w_0=1$. For any  $1\le i\le s$ we denote
$$w_i = r_{\xi_1}r_{\xi_2}\cdots r_{\xi_i},\eqno(1)$$
 The set  $\{\eta|~ \xi_i\succ \eta\}$ decomposes into a disjoint union  $$B_i\sqcup D^-_i\sqcup D_i^+.$$
{\bf Proposition 2.3}.\\ \emph{
1)~ If  $\eta\in B_i$, then  $w_i(\eta)>0$.\\
2)~ If  $\eta\in D^-_i\cup D_i^+$, then $w_i(\eta) < 0$.}
\\
{\bf Proof} will be proceeded by induction on  $0\le i\le s$. For
$i=0$ the statement is evident. Suppose that the statement is true
for all numbers that is less than  $i$. Let us prove the statement
for  $i$.

Let $\xi_i=(k,t)$, ~$k>t$.  Decompose  $\{\eta|~ \xi_i\succ \eta\}$
into four subsets  $I\sqcup II\sqcup III\sqcup IV$,
where
$$I=\{(b,t)|~  t<b<k\},$$
$$II=\{(k,c)|~ t < c < k \},$$
$$ III = \{(b,k)|~ k<b \},$$
and  $IV$ is consists of pairs  $\{\eta|~ \xi_i\succ \eta\}$, that
do not lie in  $I$, $II$ and $III$.

If $\eta\in IV$, then  $r_{\xi_i}(\eta)=\eta$, and the statement of
proposition follows from the induction assumption.
\\
{\bf Case 1}. Let $\eta\in I\cap B_i$. Then $\eta=(b,t)$, where $
t<b<k$,  and the place $\eta$ is empty after the $i$th step.
Therefore, the place  $(k,b)$ is filled  ( by the symbol  "$-$")
before the $i$th step (otherwise at the  $i$th step the place
$(b,t)$ is filled by the symbol  "$+$") . The  induction assumption
implies that  $w_{i-1}(\eta')<0$ for  $\eta'=(k,b)$. Hence,
$$w_i(\eta)=w_{i-1}r_{\xi_i}(\eta)=-w_{i-1}(\eta')>0.$$
{\bf Case 2}. Let  $\eta\in I\cap (D_i^-\cup D_i^+)$. In the
following items  2a) and  2b) we shall show that $w_i(\eta)<0$.\\
{\bf 2a)}.  Let  $\eta\in I\cap C_i^+$.  Then $\eta=(b,t)$, where $
t<b<k$,  and the place  $\eta$ is filled by the symbol  "$+$"\ at
the $i$th step. Hence, the place  $(k,b)$ is filled by the symbol
"$-$"\ at the  $i$th step and empty at the previous  $(i-1)$th step.
the induction assumption implies  $w_{i-1}(\eta')>0$ for
$\eta'=(k,b)$. Therefore,
$$w_i(\eta)=w_{i-1}r_{\xi_i}(\eta)=-w_{i-1}(\eta')< 0.$$
{\bf 2b)}. Let  $\eta\in I\cap (D_i^-\cup D_i^+)$ and $\eta\notin
C_i^+$. In this case the place $\eta$ is filled by one of the
symbols  "$+$"\ or "$-$"\ at some  $j$th step, $j<i$.

Let us show that the place  $\eta'$, that is equal to  $-
r_{\xi_i}(\eta) =(k,b)$, is free after the  $i$th step  (i.e.
$\eta'\in B_i$). Let us suppose the contrary. Then $\eta'\in C_j^-$
for some  $j<i$. The place  $\xi_j$ is situated in a column with the
number less or equal to  $t$. Since  $\eta\in D_i^-\cup D_i^+$, then
$\eta\in D_i^-$ or $\eta\in D_i^+$. From  $\eta\in D_i^-$ and
$\eta'\in C_j^-\subset D_i^-$ we have got  $\xi_i=\eta+\eta'\in
D_i^-$ (see proposition 1.4). This leads to a contradiction.

We have to treat the case  $\eta\in D_i^+$. Then  $\eta\in C_m^+$
for some  $\xi_m=(c,t)$, ~$c>k$, that situates in the same the $t$th
columns as $\xi_i$,  but below  $\xi_i$. Then at the $m$th step the
place  $\eta$, that is equal to  $(b,t)$, will be filled by the
symbol  "$+$", and the place  $(c,b)$ --- by the symbol  "$-$". At
the $m$th step the place  $(c,k)$ is already filled by the symbol
"$-$" (otherwise at the  $m$th step  $(c,k)$ will be filled by
"$-$", and $\xi_i$, that is equal to  $(k,t)$, will be filled by
"$+$"). Finally, after the  $(m-1)$th step the places $(k,b)$, that
equals to  $\eta'$, and  $(c,k)$ are already filled by the symbol
"$-$", and at the same time the place  $(c,b)$ is empty. This
contradicts to the statement that  $\dog^-_{m-1}$ is a subalgebra
(see proposition 1.4).

At last, $\eta'\in B_i$. Then  $w_{i-1}(\eta')>0$ and, therefore,
$w_i(\eta)=-w_{i-1}(\eta')<0$. This proves 2) for the case $\eta\in
I$.\\
 {\bf Case 3}.  $\eta\in II\cap B_i$.
The case is treated similarly to the case  1.\\
{\bf Case 4}.  $\eta\in II\cap (D_i^-\cup D_i^+)$. The case is
treated similarly to the case  2.\\
 {\bf Case 5}. $\eta\in III$.
As in the proof of theorem 2.2, we denote by  $a^{(t)}$ the greatest
number, that depends on  $t$, such that  $(a^{(t)},t)\notin M$. Then
all pairs  $(c,t)$, where $c>t$, lie in  $M$. Recall that  $III =
\{(b,k)|~ k<b \}$.  Decompose  $III$ into two subsets :\\
 $III_1 = \{(b,k)\in III |~ k< b \le a^{(t)}\}$,\\
  $III_2 =
\{(b,k)\in III |~ n\ge b> a^{(t)}\}$. \\
{\bf 5a}. Let  $\eta \in III_2$. By definition of  $a^{(t)}$, the
rectangle  $(a^{(t)},n]\times [1,t]$ is filled by the symbol
"$\bullet$"\ in the diagram  $\DC_\LC$. All places  $\xi_j$,~ $j\le
i$, are situated upper this rectangle. Hence  $III_2\subset B_i$ and
$w_i(\eta) >0$. This proves the statement of proposition in this
case.\\
{\bf 5b}. Let us show that  $III_1\subset \bigcup_{j<i} C_j^-$.
Suppose the contrary, let there exists  $\eta=(b,k)\in III_1$ such
that
$$ \eta \notin \bigcup_{j<i} C_j^-.$$
Then the place  $\eta$ is empty after the  $i$th step. Consider the
place  $(b,t)$. Since  $b>k$, then the place  $(b,t)$ is filled
before the  $i$th step by one of the symbols  "$\otimes$", "$-$"\ or
"$+$".

On the other hand, the symbol "$\times$"\ can not take the place
$(b,t)$, since then $\eta=\xi_m$ for some  $m<i$. At the  $m$th step
the place  $\xi_i=(k,t)$ will be filled by the symbol  "$+$"\ (resp.
 $(b,k)$ --- by the symbol  "$-$"). The place  $(b,t)$ can not be filled
 by the symbol, since in this case  $(k,t)\in B_{i-1}$, ~$(b,k)\in
B_i\subset B_{i-1}$ and $(b,t)\notin B_{i-1}$, this contradicts to
proposition  1.3 (the subspace $\nog_{i-1}$ is not a subalgebra).

Finally, suppose that the place  $(b,t)$ is filled by the symbol
"$+$". Then there exists  $\xi_j=(c,t)\in S$,~ $j<i$, ~$c>b$, such
that at the  $j$th step the symbol "$+$"\ appears on the place
$(b,t)$   and the symbol "$-$"\  on the place  $(c,b)$. At the same
 $j$th step the place  $(c,k)$ must be already filled, otherwise  $(k,t)$
would be filled by the symbol  "$+$"\ after the $i$th step.
According to the procedure of arrangement  of symbols,  the place
$(c,k)$ may be filled only by the symbol  "$-$". Thus, at the $j$th
step we have got that the places  $(c,b)$ and  $(b,k)$ are empty,
and the place $(c,k)$ is filled by the symbol  "$-$". this
contradicts to the
statement that  $\nog_{j-1}$ is a subalgebra (see proposition 1.4). The statement of the item  5b is proved. \\
{\bf 5c}. Let us show that $w_i(\eta)<0$ for any  $\eta=(b,k)\in
III_1$. Let  $\eta=(b,k)$. Then  $(b,t)=r_{\xi_i}(b,k)$. Let  $m$ be
the greatest number such that  $\xi_m\succeq (b,t)$.  The following
two cases is possible.\\
 {\bf a)}~ $m=i-1$. Then either  $(b,t)=\xi_{i-1}$, or $(b,t)\in
 D_{i-1}^+ \sqcup D_{i-1}^-$. In any case  $w_{i-1}(b,t)<0$. We have got
 $$ w_i(\eta)=w_{i-1}r_{\xi_{i}}(b,k))= w_{i-1}(b,t))<0.$$
 {\bf b)}~ $m > i-1$.  Let  $\xi_{i-1} = (c,t)$. Then
$b>c$ and  $$r_{\xi_{i-1}}r_{\xi_{i}}(b,k)=
 r_{\xi_{i-1}}(b,t)=(b,c).$$
Note that  $r_{\xi_p}(b,c)=(b,c)$ for all  $ m<p<i-1$. Hence,
$$ w_i(\eta)= w_mr_{\xi_{m+1}}\ldots r_{\xi_{i-1}}r_{\xi_{i}}(b,k) =
w_m(b,c).\eqno(2)$$ We have to prove that $w_m(b,c)<0$.

Since in the pair  $\xi_m$ the number of column is less or equal to
 $t$, and  $c>t$, then after the $m$th step
 the place  $(b,c)$ is either empty, or filled by the symbol "$-$". By induction assumption, the inequality  $w_m(b,c)<0$
 is equivalent to the statement that the place
$(b,c)$ is filled after the  $m$th step (by the symbol  "$-$").

Suppose that the place  $(b,c)$ is empty after the  $m$th step.
 The place  $(b,t)$ is filled after the  $i$th step, as long as  $b>k$. By definition of  $m$,
the place $(b,t)$ is filled after the  $m$th step (by one pf the
symbols  $\otimes$, ~"$+$", ~"$-$"). So, after the  $m$th step the
place  $(b,t)$ is already filled, and the places  $(b,c)$ and
$(c,t)$ are empty. This contradicts to the statement that $\nog_m$
is a subalgebra. Thus, $(b,c)$ is already filled after the  $m$th
step by the symbol  "$-$"\ and, therefore, $w_m(b,c)<0$. This proves
that  $w_i(\eta)<0$.
$\Box$\\
{\bf Corollary 2.4}. \emph{If $\eta\in B_i$, then $w_j(\eta)>0$ for
any
$1\le j\le i$.}\\
{\bf Proof} follows from the inclusion  $B_i\subset B_j$. $\Box$

Denote by  $w^{(t)}$ the product  of reflections $r_\xi$ (arranged
in the  decreasing order in the sense of $\succ$), where the number
of column of $\xi$ is equal to $t$.
 We construct the system of permutations
$$w^{[t]}= w^{(1)}\ldots w^{(t)}.\eqno(2)$$ Remark that
$w^{[t]}$  coincides with  $w_i$, where $\xi_i$ is the least (in the
sense of $\succ$) root among all roots lying in the first  $t$
columns.

Let as above $a^{(t)}=\max\{c|~ (c,t)\notin M\}$.\\
 {\bf Proposition 2.5}. \emph{We claim that\\
1)~ $w^{[t]}(\eta) > 0$ for any  $\eta=(b,t)$, ~ $ a^{(t)}< b\le n$.\\
2)~  $w^{[t]}(\eta) < 0$ for any $\eta=(b,t)$,
 ~ $t<b\le a^{(t)}$.}\\
 {\bf Proof}.
By definition of $a^{(t)}$, the rectangle   $(a^{(t)},n]\times
[1,t]$ is filled by "$\bullet$"\ in the diagram $\DC_\LC$. All
places
 $\xi_j$,~ $j\le i$, are situated upper this rectangle. This implies
 the statement  1).

 We shall prove  the statement  2) in each of these cases separately.
 \\
{\bf i)} The $t$th column does not contain the root of $S$. In this
case all column is filled by the symbol  "$-$"\.~ We have got
$w^{[t]}(\eta)=w^{[t-1]}(\eta)$. By proposition  2.3, we obtain
$w^{[t-1]}(\eta)<0$.\\
{\bf ii)} Let   $(b,t)$ situate upper all roots from  $S$, lying in
the  $t$th column. Then  $(b,t)$ is filled by the symbol  "$+$"\ or "$-$". By proposition 2.3, $w^{[t]}(\eta) < 0$.\\
{\bf iii)} $(b,t)$ situates lower  $\xi_i$, that is the least in the
sense of  $\succ$ root of the  $t$th column, or coincides with it.

If  $\eta=\xi_i$, then
$$w^{[t]}(\eta)=w_i(\xi_i)=-w_{i-1}(\xi_i)<0.$$

Let  $\xi_i=(k,t)$ and $b>k$. Then  $r_{\xi_i}(b,t)=(b,k)$. Let
 $\xi_m$  be the least root of  $S$, that is greater (in the sense of  $\succ$) or equal to  $(b,t)$.
By item  5b (see the proof of proposition  2.3),  $(b,k)$ is filled
by the symbol  "$-$"\  before or during the  $m$th step. Therefore,
~$w_{m}(b,k)<0$. Since  $r_{\xi_p}(b,k)=(b,k)$ for any  $m<p<i$,
then
$$w^{[t]}(\eta) = w_i(\eta)=w_{m}r_{\xi_{m+1}}\ldots r_{\xi_{i}}(b,t)=
w_{m}r_{\xi_{m+1}}\ldots r_{\xi_{i-1}}(b,k)=  w_{m}(b,k)  < 0.$$
 $\Box$

 Let as above  $S=\xi_1\succ\ldots\succ \xi_s\}$.
 Recall that  $(k,t)\in S$ if and only if the place
 $(k,t)$ is filled by the symbol "$\otimes$"{} in the diagram  $\DC_\LC$.\\
 {\bf
Theorem  2.6}. \emph{We claim that   $w = r_{\xi_1}r_{\xi_2}\cdots
r_{\xi_s}$.}\\
{\bf Proof}. according to formula  (2) $$w^{[n]} =
r_{\xi_1}r_{\xi_2}\cdots r_{\xi_s}.$$

 Let us prove  $w^{[n]}(t) =w(t)$ for any $1\le t\le n$, using induction on $t$.
 For  $t=1$ the statement is evident. Suppose that the statement is true for the numbers less than  $t$. Now we shall prove it for
  $t$.

By definition, $w(t)$ is a greatest number among the numbers of the
segment  $[1,n]$ of positive integers that do not lie in

$$\Lambda_t = \{w(1), \ldots, w(t-1)\}\sqcup (a^{[t]}, n].\eqno(3)$$

By induction assumption, $w(j)=w^{[n]}(j)$ for any  $1\le j\le t-1$.
Note that  $w^{[n]}(j)=w^{[t]}(j)$, as long as  $r_\xi(j)=j$ for all
 $\xi\in S$, lying in the columns with numbers greater than  $t$.
Hence, $w(j)= w^{[t]}(j)$ for all $1\le j\le t-1$. This gives a
possibility to substitute $w$  for  $w^{[t]}$ in (3).

By definition of  $a^{(t)}$, the rectangle  $(a^{(t)},n]\times
[1,t]$ is filled by the symbol "$\bullet$"\ in the diagram
$\DC_\LC$. The places  $\xi_j$,~ $j\le i$ situates upper this
rectangle. Hence, for any $p\in (a^{[t]}, n]$ we have
$w^{[t]}(p)=p$.

 The set  $\Lambda_t$ may be represented in the form
 $$
\Lambda_t = \{w^{[t]}(1), \ldots,
w^{[t]}(t-1)\}\sqcup\{w^{[t]}(p)\vert ~ a^{[t]}< p\le n\}.\eqno(4)$$
All elements of the segment  $[1,n]$ of positive integers, that do
not lie in  $\Lambda_t$, have the form  $w^{[t]}(k)$, ~ where $t\le
k\le a^{(t)}$.

From item  2) of proposition  2.5, ~ $w^{[t]}(t)>w^{[t]}(k)$ for any
 $t<k\le a^{(t)}$. Therefore, $w^{[t]}(t)$ is a greatest number among all numbers of the  segment
  $[1,n]$ of positive integers, that do not lie in  $\Lambda_t$. We conclude that  $w^{[t]}(t)=w(t)$. Finally,
$w^{[t]}(t) =w^{[n]}(t)$, as long as  $r_\xi(t)=t$ for all  $\xi\in
S$, lying in columns with numbers greater than  $t$. At last,
$w^{[n]}(t) =w(t)$. $\Box$
\\
Denote by  $A^{(t)}$ the set   $\eta\in A$, that has the form
$(b,t)$ for some  $b>t$. \\
{\bf Theorem  2.7}. \emph{Let  $\eta\in A^{(t)}$, then \\
1)~ the place  $\eta$ is filled in the diagram  $\DC_\LC$ by the
symbol  "$-$"\ iff
 $w^{[t-1]}(\eta)<0$;\\
2)~ the place  $\eta$ is filled in the diagram  $\DC_\LC$ by the
symbol
"$\bullet$"{} iff  $w^{[t]}(\eta)>0$;\\
3)~ место $\eta$ is filled in the diagram  $\DC_\LC$ by the symbol
"$+$"\ or "$\otimes$"{} iff $w^{[t-1]}(\eta)>0$ и $w^{[t]}(\eta)<
0.$}\\
 {\bf Proof } is a corollary of the propositions 2.3 and 2.5. $\Box$

\section*{3. Field of invariants }

In this section we formulate the conjecture on the structure of
field of invariants of the coadjoint representation of the Lie
algebra $\LC$.

To any $\xi\in S$ we shall correspond a polynomial  $P_\xi$.
 Let   $\xi =\xi_m = (k,t)\in S$, where
$k>t$. Denote by  $w_\xi$ the permutation  $w_m =r_{\xi_1}\ldots
r_{\xi_m}$.

 {\bf Case  1}. $w_\xi(t)>t$. On can show that in this case  $w_\xi(t)=k$. Put
$$J:=J(\xi)=\{1\le j\le t:~w_\xi(j)\ge w_\xi(t)\}, \quad\quad I:=I(\xi) =
wJ(\xi) .$$
 {\bf Case  2}. $w_\xi(t)\le t$. The system
$J:=J(\xi)$ is defined as in  (3.1). Denote
$$ I_*(\xi):=I_*(\xi)=\{1\le i\le n:~ i>t, ~w_\xi(i) < w_\xi(t)\},$$
$$I(\xi)=[w_\xi(t),t]\sqcup I_*(\xi).$$
On can show that in both cases  $|I(\xi)|=|J(\xi)|$.

As above  $\{y_{ij}\}$ is the standard basis in $\nog$. By the
diagram  $\DC_\LC$, we construct the matrix  $\Phi_\LC$, in which
the places  $\{(i,j)\in \Dp\setminus M\} $ are filled by the
corresponding elements  $\{y_{ij}\}$ of the standard basis; the
other places are filled by zeroes. For instance, for the Lie algebra
 $\LC$ of the example 1, we have the following diagram   $\DC_\LC$ and matrix
$\Phi_\LC$:

\begin{tabular}{cc}{\large $\DC_\LC$ =
\begin{tabular}
{|p{0.1cm}|p{0.1cm}|p{0.1cm}|p{0.1cm}|p{0.1cm}|p{0.1cm}|p{0.1cm}|}
\hline  &  &  &  & & &  \\
\hline $+$& & &  & & &  \\
\hline $+$&$+$ & & & & &  \\
\hline $\otimes$ & $-$ & $-$  &  & &  & \\
\hline $\bullet$ & $+$ & $+$&$\otimes$ & & &\\
\hline  $\bullet$ & $\otimes$ & $-$ &$+$  &$-$ & &\\
\hline   $\bullet$&  $\bullet$ & $\otimes$ &$\otimes$  &$-$  &$-$ & \\
\hline
\end{tabular}}\quad,&

$ \quad\quad
 \Phi_\LC=\left(\begin{array}{ccccccc}
0&0&0&0&0&0&0\\
y_{21}&0&0&0&0&0&0\\
 y_{31}&y_{32}&0&0&0&0&0\\
  y_{41}&y_{42}&y_{43}&0&0&0&0\\
   0&y_{52}&y_{53}&y_{54}&0&0&0\\
    0&y_{62}&y_{63}&y_{64}&y_{65}&0&0\\
     0&0&y_{73}&y_{74}&y_{75}&y_{76}&0
\end{array}\right)
$
\end{tabular}

Let  $\la$  be a variable. Consider the characteristic matrix
$\Phi_\LC - \la E$. Note that any minor of the characteristic matrix
is a polynomial in  $\la$ with coefficients in $S(\LC)=K[\LC^*]$.

Let   $M_\xi(\la)$ be the polynomial of characteristic matrix with
the system of columns  $J(\xi)$ and the system of row  $I(\xi)$, and
 $P_\xi$ be its highest coefficient.\\
{\bf Conjecture}. The field of invariants of the coadjoint
representation of the Lie algebra $\LC$ is a field of rational
functions in  $P_\xi$,~ $\xi\in S$.

 In the example  1, $S$ consists of four elements
 $\{\xi_1, \xi_2, \xi_3, \xi_4\}$.
Direct calculations shows that  $P_{\xi_1} = y_{41}$,~ $ P_{\xi_2} =
y_{62}$,~ $P_{\xi_3} = y_{73}$,~ $ P_{\xi_4} = y_{74}y_{41}
+y_{73}y_{31}$.

Samara State University\\
443011, Samara\\
ul. akad. Pavlova, 1,\\
Russia\\
\emph{ E-mail}: apanov@list.ru

\end{document}